\documentclass[11pt]{amsart}

\usepackage{amsmath}
\usepackage{amssymb}
\usepackage{mathrsfs}

\newtheorem{theorem}{Theorem}[section]
\newtheorem{claim}[theorem]{Claim}

\newtheorem{lemma}[theorem]{Lemma}

\newtheorem{corollary}[theorem]{Corollary}

\theoremstyle{definition}
\newtheorem{definition}[theorem]{Definition}

\newtheorem{question}[theorem]{Question}

\theoremstyle{remark}

\newcount\skewfactor
\def\mathunderaccent#1#2 {\let\theaccent#1\skewfactor#2
\mathpalette\putaccentunder}
\def\putaccentunder#1#2{\oalign{$#1#2$\crcr\hidewidth
\vbox to.2ex{\hbox{$#1\skew\skewfactor\theaccent{}$}\vss}\hidewidth}}

\def\smallbox#1{\leavevmode\thinspace\hbox{\vrule\vtop{\vbox
   {\hrule\kern1pt\hbox{\vphantom{\tt/}\thinspace{\tt#1}\thinspace}}
   \kern1pt\hrule}\vrule}\thinspace}

\newcommand{\cf}{{\rm cf}}

\def\qedref#1{$\qed_{\reforiginal{#1}}$}

\title{Dear lambda}
\author{Shimon Garti}
\address{Institute of Mathematics
 The Hebrew University of Jerusalem,
 Jerusalem 91904, Israel}
\email{shimon.garty@mail.huji.ac.il}

\author{Saharon Shelah}
\address{Institute of Mathematics
 The Hebrew University of Jerusalem,
 Jerusalem 91904, Israel
 and  Department of Mathematics
 Rutgers University
 New Brunswick, NJ 08854, USA}
\email{shelah@math.huji.ac.il}
\urladdr{http://www.math.rutgers.edu/\char`\~shelah}
\thanks{Research supported by European Research Council, grant 338821. This is publication 1143 of the second author}


\begin{document}
\let\labeloriginal\label
\let\reforiginal\ref
\def\ref#1{\reforiginal{#1}}
\def\label#1{\labeloriginal{#1}}

\begin{abstract}
We prove the consistency of $\mathfrak{r}_\lambda<\mathfrak{d}_\lambda$, and even $\mathfrak{u}_\lambda<\mathfrak{d}_\lambda$, for a singular cardinal $\lambda$.
\end{abstract}

\maketitle

\newpage

\section{Introduction}

This paper deals with several cardinal characteristics of the continuum, including the reaping number and the ultrafilter number.
We define them in the generalized form of $\mathfrak{r}_\kappa$ and $\mathfrak{u}_\kappa$ where $\kappa$ is any infinite cardinal.

\begin{definition}
\label{defrkappa} The reaping number. \newline 
Let $\kappa$ be an infinite cardinal, $B\in[\kappa]^\kappa$.
\begin{enumerate}
\item [$(\aleph)$] A set $S\in[\kappa]^\kappa$ splits $B$ iff $|S\cap B|=|(\kappa-S)\cap B|=\kappa$.
\item [$(\beth)$] A family of sets $\mathcal{A}\subseteq[\kappa]^\kappa$ is called an unreaped family iff there is no single $S\in[\kappa]^\kappa$ which splits any element of $\mathcal{A}$.
\item [$(\gimel)$] The reaping number $\mathfrak{r}_\kappa$ is the minimal cardinality of an unreaped family in $[\kappa]^\kappa$.
\end{enumerate}
\end{definition}

An unreaped family will also be called an unsplittable or an $\mathfrak{r}_\kappa$ family.
A closed friend of the reaping number is the ultrafilter number. 
Recall that an ultrafilter $\mathscr{U}$ over $\kappa$ is uniform iff all of its elements are of size $\kappa$, and all the ultrafilters in this paper are uniform.
In particular, the definition of the ultrafilter number applies to uniform ultrafilters.
Again, we phrase the definition in the general context.

\begin{definition}
\label{defukappa} The ultrafilter number. \newline 
Let $\kappa$ be an infinite cardinal and let $\mathscr{U}$ be an ultrafilter over $\kappa$.
\begin{enumerate}
\item [$(\aleph)$] A base for $\mathscr{U}$ is a collection $\mathcal{B}\subseteq\mathscr{U}$ such that for any $A\in\mathscr{U}$ there is some $B\in\mathcal{B}$ so that $B\subseteq A$.
\item [$(\beth)$] The characteristic of a uniform ultrafilter $\mathscr{U}$ over $\kappa$ is the minimal size of a base for $\mathscr{U}$, denoted by ${\rm Ch}(\mathscr{U})$.
\item [$(\gimel)$] The ultrafilter number $\mathfrak{u}_\kappa$ is the minimal size of a base for some uniform ultrafilter $\mathscr{U}$ over $\kappa$.
\end{enumerate}
\end{definition}

Lest $\kappa=\aleph_0$ we denote $\mathfrak{r}_\kappa$ by $\mathfrak{r}$ and $\mathfrak{u}_\kappa$ by $\mathfrak{u}$.
Any base of an ultrafilter is unsplittable, hence $\mathfrak{r}_\kappa\leq\mathfrak{u}_\kappa$.
It is easy to see that both $\mathfrak{r}_\kappa>\kappa$ and $\mathfrak{u}_\kappa>\kappa$.
Our purpose in the first section is to analyze the cofinality of these characteristics.

Cardinal characteristics which may assume countable cofinality are rarefied, and for proving this fact one needs complicated arguments. 
The almost disjointness number $\mathfrak{a}$ is an example, as proved in \cite{MR1975392}.
It is unknown whether the cofinality of $\mathfrak{r}_\kappa$ is always greater than $\kappa$, and in particular whether $\cf(\mathfrak{r})$ is uncountable.
Additional examples are $\mathfrak{i}$ (see \cite{MR1975392}) and $\mathfrak{gp}$ (see \cite{MR3787522}, the possibility of countable cofinality in this case requires instances of Chang's conjecture hence depends on the existence of large cardinals).
The question about $\mathfrak{r}$ appeared in \cite{MR1234292}, Problem 3.4.
By adding $\lambda$ many Cohen reals to a model of GCH we obtain $\mathfrak{r}=\mathfrak{u}=\mathfrak{c}=\lambda$ and hence these characteristics may be singular.
We shall prove that the cofinality of $\mathfrak{u}_\kappa$ is above $\omega$, and in some cases the same holds for $\mathfrak{r}_\kappa$.
Namely, if $\kappa=\cf(\kappa)$ and $\mathfrak{r}_\kappa<\mathfrak{d}_\kappa$ then $\mathfrak{r}_\kappa=\mathfrak{u}_\kappa$ and hence $\cf(\mathfrak{r}_\kappa)>\omega$.

This brings us to the question whether $\mathfrak{r}_\kappa<\mathfrak{d}_\kappa$ is possible.
If $\kappa=\aleph_0$ then the answer is positive.
By adding $\omega_2$ Miller reals to a model of the continuum hypothesis one obtains $\mathfrak{r}=\omega_1<\omega_2=\mathfrak{d}$.
Raghavan and Shelah, \cite{RagSh}, proved that if $\kappa=\cf(\kappa)>\beth_\omega$ then $\mathfrak{d}_\kappa\leq\mathfrak{r}_\kappa$.
We shall prove in the second section that $\mathfrak{r}_\lambda<\mathfrak{d}_\lambda$ is consistent for a strong limit singular cardinal $\lambda$.
This can be done with any cofinality of $\lambda$, and above $\beth_\omega$.

Our notation is standard.
We mention here the concept of strong finite intersection property.
Let $\kappa$ be an infinite cardinal, $\mathscr{F}\subseteq[\kappa]^\kappa$.
We say that $\mathscr{F}$ has the strong finite intersection property iff $|\bigcap u|=\kappa$ whenever $u\in[\mathscr{F}]^{<\omega}$.
We suggest \cite{MR2768685} as an excellent background regarding cardinal characteristics.
For a good background in pcf theory we suggest \cite{MR2768693}.

\newpage 

\section{Cofinality}

We open this section with a theorem about the cofinality of the ultrafilter number. The statement and the proof are phrased in the case of $\kappa=\aleph_0$, and possible generalizations to higher cardinals are discussed after the proof.

\begin{theorem}
\label{thmmmt} The cofinality of a base. \newline 
Let $\mathscr{U}$ be a non-principal ultrafilter over $\omega$. \newline 
If ${\rm Ch}(\mathscr{U})=\mu$ then $\cf(\mu)>\omega$. \newline 
\end{theorem}

\par\noindent\emph{Proof}. \newline 
Fix an ultrafilter $\mathscr{U}$ over $\omega$ such that ${\rm Ch}(\mathscr{U})=\mu$.
Assume by way of contradiction that $\cf(\mu)=\omega$.
Choose an increasing sequence of uncountable ordinals $(\mu_n:n\in\omega)$ such that $\mu=\bigcup_{n\in\omega}\mu_n$.
Fix a base $\mathcal{B}=\{B_\beta:\beta\in\mu\}$ for the ultrafilter $\mathscr{U}$.

For every $\alpha\in\mu$ let $\mathcal{B}_\alpha=\{B_\beta:\beta<\alpha\}$.
By our assumption toward contradiction we see that each $\mathcal{B}_\alpha$ fails to be a base for $\mathscr{U}$.
Consequently, for every $\alpha\in\mu$ one can choose $y_\alpha\in[\omega]^\omega-\mathscr{U}$ such that:
\begin{enumerate}
\item [$(\alpha)$] $u\in[\alpha]^{<\omega} \Rightarrow |\bigcap_{\beta\in u}B_\beta\cap y_\alpha|=\aleph_0$.
\item [$(\beta)$] $u\in[\alpha]^{<\omega} \Rightarrow |y_\alpha-\bigcup_{\beta\in u}y_\beta|=\aleph_0$.
\end{enumerate}
Indeed, for every $\alpha\in\mu$ since $\mathcal{B}_\alpha$ does not generate $\mathscr{U}$ there is some $x_\alpha\in\mathscr{U}$ such that for every $u\in[\alpha]^{<\omega}$ it is true that $\neg(\bigcap_{\beta\in u}B_\beta\subseteq^* x_\alpha)$.
Let $y_\alpha=\omega-x_\alpha$ and conclude that $(\alpha)$ holds. Item $(\beta)$ is just an equivalent formulation of the same statement.

Using the above property we define, by induction on $n\in\omega$, a set $s_n$ for which the following requirements are met:
\begin{enumerate}
\item [$(a)$] $s_n$ is an infinite subset of $\omega$.
\item [$(b)$] $s_n\notin\mathscr{U}$.
\item [$(c)$] $s_n$ is disjoint from $\bigcup_{m<n}s_m$.
\item [$(d)$] $n\in\bigcup_{m\leq n}s_m$.
\item [$(e)$] $u\in[\mu_n]^{<\omega}\Rightarrow |\bigcap_{\beta\in u}B_\beta\cap s_n|=\aleph_0$.
\end{enumerate}
The choice can be done, basically, by $(\alpha)$ (or $(\beta)$) above.
If $n=0$ then let $\alpha=\mu_0$ and choose $y_\alpha$ as guaranteed in $(\alpha)$ and $(\beta)$. 
Set $s_0 = y_\alpha\cup\{0\}$.
Observe that $(a),(b)$ are satisfied and $(c)$ is vacuous in this case. We added zero to $y_\alpha$ in order to satisfy $(d)$, and $(e)$ is exactly $(\alpha)$ with respect to $y_\alpha$ hence also to $s_0$.

In the $(n+1)$st stage we choose a sufficiently large $\ell>n$ and let $\alpha=\mu_\ell$.
Again, let $y_\alpha$ be as guaranteed in $(\alpha)$ and $(\beta)$.
Define $t_\alpha=y_\alpha-\bigcup_{m\leq n}s_m$.
From $(\beta)$ we infer that $t_\alpha$ is infinite.
Now if $n+1\in\bigcup_{m\leq n}s_m$ then let $s_{n+1}=t_\alpha$ and if $n+1\notin\bigcup_{m\leq n}s_m$ then let $s_{n+1}=t_\alpha\cup\{n+1\}$.
One can verify that all the requirements are satisfied.

For every $i\in\{0,1\}$ define $E_i = \bigcup\{s_{2n+i}:n\in\omega\}$.
From $(d)$ one can see that $E_0\cup E_1=\omega$, and from $(c)$ it follows that $E_0\cap E_1=\varnothing$, so $\{E_0,E_1\}$ is a partition of $\omega$.
Consequently, there must be some $i\in\{0,1\}$ for which $E_i=\varnothing\ \text{mod}\ \mathscr{U}$.

Let $u\subseteq\mu$ be any finite set of ordinals.
Pick up a sufficiently large $n$ so that $u\subseteq\mu_n<\mu_{2n}$.
Apply $(e)$ and conclude that:
$$
|\bigcap_{\beta\in u}B_\beta\cap s_{2n}|=|\bigcap_{\beta\in u}B_\beta\cap s_{2n+1}|=\aleph_0.
$$
By the definition of $E_0$ and $E_1$ we infer that for every $u\in[\mu]^{<\omega}$ it is true that $|\bigcap_{\beta\in u}B_\beta\cap E_0|=|\bigcap_{\beta\in u}B_\beta\cap E_1|=\aleph_0$.
But this means that $E_0\neq\varnothing\ \text{mod}\ \mathscr{U} \wedge E_1\neq\varnothing\ \text{mod}\ \mathscr{U}$, a contradiction.

\hfill \qedref{thmmmt}

The above theorem can be generalized to higher cardinals in the following manner.
Assume that $\kappa$ is an infinite cardinal and $\mathscr{U}$ is a uniform ultrafilter over $\kappa$.
If ${\rm Ch}(\mathscr{U})=\mu$ then $\cf(\mu)>\omega$, by the same proof.
We conclude, therefore, that $\cf(\mathfrak{u}_\kappa)>\omega$ at every infinite 
cardinal $\kappa$.
A stronger generalization requires some degree of completeness.
Assuming that $\kappa$ is measurable and $\mathscr{U}$ is $\kappa$-complete, the above arguments show that if ${\rm Ch}(\mathscr{U})=\mu$ then $\cf(\mu)>\kappa$.

\begin{question}
\label{qcofu} Is it consistent that $\cf(\mathfrak{u}_\kappa)\leq\kappa$ for some infinite cardinal $\kappa$?
\end{question}

Back to Theorem \ref{thmmmt}, we know that the ultrafilter number (at any cardinal) has uncountable cofinality.
Can we prove a similar theorem about $\mathfrak{r}$?
It has been shown by Aubrey, \cite{MR2058185} that if $\mathfrak{r}<\mathfrak{d}$ then $\mathfrak{r}=\mathfrak{u}$.
Hence in this case, the cofinality of the reaping number will be uncountable.
Our next goal is to generalize this result to regular uncountable cardinals.
We follow in the footsteps of Aubrey, with the required adaptations to the general case.
Let us begin with another cardinal characteristic:

\begin{definition}
\label{defdkappa} The dominating number. \newline 
Let $\kappa$ be an infinite cardinal. 
\begin{enumerate}
\item [$(\aleph)$] For $f,g\in{}^\kappa\kappa$ we shall say that $g$ dominates $f$ iff $\{\beta\in\kappa:f(\beta)>g(\beta)\}$ is of size less than $\kappa$. This relation will be denoted by $f\leq^* g$.
\item [$(\beth)$] A family of functions $\mathcal{D}\subseteq{}^\kappa\kappa$ is called a dominating family iff for every $f\in{}^\kappa\kappa$ there exists $g\in\mathcal{D}$ such that $f\leq^* g$.
\item [$(\gimel)$] The dominating number $\mathfrak{d}_\kappa$ is the minimal size of a dominating family at ${}^\kappa\kappa$.
\end{enumerate}
\end{definition}

Let $\kappa$ be a regular cardinal. We shall say that $\Pi = \{I^\Pi_\alpha: \alpha<\kappa\}$ is an interval partition of $\kappa$ when each $I^\Pi_\alpha$ is a non-empty interval of the form $[\gamma_\alpha,\gamma_{\alpha+1})$, if $\alpha<\beta<\kappa$ then $I^\Pi_\alpha\cap I^\Pi_\beta=\varnothing$ and every ordinal of $\kappa$ belongs to some $I^\Pi_\alpha$.
If $\Pi$ is clear from the context then we may write $I_\alpha$ instead of $I^\Pi_\alpha$.

\begin{definition}
\label{defintervals} Let $\Pi=\{I_\alpha:\alpha<\kappa\}$ be an inverval partition of $\kappa$. Assume that $\mathscr{F}\subseteq[\kappa]^\kappa$.
\begin{enumerate}
\item [$(\aleph)$] A pair $(D,E)$ is a nice $\Pi$-orbit iff both $D$ and $E$ are unions of $\kappa$-many intervals from $\Pi$, the intervals of $D$ are disjoint from the intervals of $E$ and moreover there is no interval of $D$ which has an adjacent interval of $\Pi$ in $E$.
\item [$(\beth)$] We say that $\mathscr{F}$ is $\Pi$-scattered iff for every nice $\Pi$-orbit $(D,E)$ and every $y\in\mathscr{F}$, both $y\cap D$ and $y\cap E$ are of size $\kappa$. 
\item [$(\gimel)$] We say that $\Pi$ is $\mathscr{F}$-scattered iff one can find a nice $\Pi$-orbit $(D,E)$ such that $y\cap D$ and $y\cap E$ are of size $\kappa$ for every $y\in\mathscr{F}$.
\end{enumerate}
\end{definition}

The following gives a simple example.
We use the terminology of \emph{almost every} in the sense that the set of exceptions is of size less than $\kappa$.

\begin{lemma}
\label{lemcase1} Let $\Pi=\{I_\alpha:\alpha<\kappa\}$ be an inverval partition of $\kappa$, and assume that $\mathscr{F}\subseteq[\kappa]^\kappa$. \newline 
If every $y\in\mathscr{F}$ meets almost every interval of $\Pi$ then $\mathscr{F}$ is $\Pi$-scattered.
\end{lemma}

\hfill \qedref{lemcase1}

Suppose that $\Pi$ is an interval partition and $\mathscr{F},\mathscr{G}\subseteq[\kappa]^\kappa$.
In the theorem below it is shown that if $\mathscr{F}$ is not $\Pi$-scattered and $\Pi$ is not $\mathscr{G}$-scattered, then one can define an interesting $\mathfrak{r}_\kappa$ family out of $\mathscr{F}$ and $\mathscr{G}$.

\begin{theorem}
\label{thmmt} Assume that:
\begin{enumerate}
\item [$(\aleph)$] $\mathscr{F},\mathscr{G}\subseteq[\kappa]^\kappa$.
\item [$(\beth)$] $\Pi = \{I_\alpha:\alpha<\kappa\}$ is an interval partition of $\kappa$.
\end{enumerate}
Then at least one of the following obtains:
\begin{enumerate}
\item [$(a)$] $\mathscr{F}$ is $\Pi$-scattered.
\item [$(b)$] $\Pi$ is $\mathscr{G}$-scattered.
\item [$(c)$] There exist $y\in\mathscr{F}$ and $h\in{}^\kappa\kappa$ which is increasing and $<\kappa$-to-one such that $\{h(y\cap z):z\in\mathscr{G}\}$ is an $\mathfrak{r}_\kappa$-family.
\end{enumerate}
\end{theorem}

\par\noindent\emph{Proof}. \newline 
If every $y\in\mathscr{F}$ meets almost every interval of $\Pi$ then $\mathscr{F}$ is $\Pi$-scattered by Lemma \ref{lemcase1}, so we may assume that this is not the case and fix some $y\in\mathscr{F}$ which evades $\kappa$-many intervals of $\Pi$.
We create a new interval partition $\Phi = \Phi(\Pi)$ by defining the intervals $I^\Phi_\alpha$ using induction on $\alpha\in\kappa$.
Every interval $I^\Phi_\alpha$ will be a union of intervals from $\Pi$.
This will be a union of consecutive intervals with a last element.

If $\alpha=0$ then let $\gamma_0\in\kappa$ be the first ordinal for which $y\cap I^\Pi_{\gamma_0}=\varnothing$ and let $I^\Phi_0 = \bigcup\{I^\Pi_\beta: \beta\leq\gamma_0\}$.
In the stage of $\alpha+1$ we assume that $I^\Phi_\alpha$ is at hand and let $I^\Pi_{\gamma_\alpha}$ be the last interval from $\Pi$ in $I^\Phi_\alpha$.
Let $\gamma_{\alpha+1}\in\kappa$ be the first ordinal greater than $\gamma_\alpha$ such that $y\cap I^\Pi_{\gamma_{\alpha+1}}=\varnothing$, and define $I^\Phi_{\alpha+1} = \bigcup\{I^\Pi_\beta: \gamma_\alpha<\beta\leq\gamma_{\alpha+1}\}$.
Notice that $I^\Phi_{\alpha+1}$ has a last interval from $\Pi$.
Finally, assume that $\alpha$ is a limit ordinal and let $\gamma = \bigcup_{\beta<\alpha}\gamma_\beta<\kappa$.
Let $\gamma_\alpha\in\kappa$ be the first ordinal bigger than $\gamma$ such that $y\cap I^\Pi_{\gamma_\alpha}=\varnothing$ and let $I^\Phi_\alpha = \bigcup\{I^\Pi_\beta: \gamma\leq\beta\leq\gamma_\alpha\}$.

We produce from the interval partition $\Phi$ a function $h:\kappa\rightarrow\kappa$ by letting $h$ be constant over the intervals of $\Phi$. Formally, $h(\delta)=\alpha$ iff $\delta\in I^\Phi_\alpha$ for every $\delta\in\kappa$.
Observe that $h$ is increasing and $<\kappa$-to-one.
Equipped with $h$ we split the rest of the proof into two cases. \newline 

\par\noindent\emph{Case 1}:
There exists $Q\subseteq\kappa$ such that $|Q|=|\kappa-Q|=\kappa$, and for every $z\in\mathscr{G}$ it is true that $|y\cap z\cap h^{-1}(Q)|=|y\cap z\cap h^{-1}(\kappa-Q)|=\kappa$. \newline 

In this case, $\Pi$ is $\mathscr{G}$-scattered so $(b)$ holds.
For proving this statement let $D = h^{-1}(Q)-\bigcup_{\alpha<\kappa}I^\Pi_{\gamma_\alpha}$ and let $E = h^{-1}(\kappa-Q)-\bigcup_{\alpha<\kappa}I^\Pi_{\gamma_\alpha}$.
Notice that $D\cap E=\varnothing$ and there are no adjacent intervals in $D,E$ (for obtaining this goal we moved out the intervals $I^\Pi_{\gamma_\alpha}$).
Fix any $z\in\mathscr{G}$. One can see that:
$$
D\cap z = (h^{-1}(Q)-\bigcup_{\alpha<\kappa}I^\Pi_{\gamma_\alpha})\cap z\supseteq (h^{-1}(Q)-\bigcup_{\alpha<\kappa}I^\Pi_{\gamma_\alpha})\cap(z\cap y).
$$
But $\bigcup_{\alpha<\kappa}I^\Pi_{\gamma_\alpha}\cap(z\cap y)=\varnothing$ by the choice of $y$, and hence $D\cap z = h^{-1}(Q)\cap z\cap y$.
By the assumption of the present case, the size of $D\cap z$ is $\kappa$.
An identical argument shows that $E\cap z$ is of size $\kappa$, upon replacing $Q$ by $\kappa-Q$.
This shows that $\Pi$ is $\mathscr{G}$-scattered as claimed in $(b)$. \newline 

\par\noindent\emph{Case 2}:
For every $Q\subseteq\kappa$ such that $|Q|=|\kappa-Q|=\kappa$, there exists some $z\in\mathscr{G}$ such that either $|y\cap z\cap h^{-1}(Q)|<\kappa$ or $|y\cap z\cap h^{-1}(\kappa-Q)|<\kappa$. \newline 

In this case we will try to create an unsplittable family out of $\mathscr{F}$ and $\mathscr{G}$ thus proving $(c)$.
Fix any $Q\subseteq\kappa$ and a set $z\in\mathscr{G}$ whose existence is guaranteed by the assumption of this case.
If $|y\cap z\cap h^{-1}(Q)|<\kappa$ then $|Q\cap h(y\cap z)|<\kappa$ and if $|y\cap z\cap h^{-1}(\kappa-Q)|<\kappa$ then $|(\kappa-Q)\cap h(y\cap z)|<\kappa$. In any case, $Q$ fails to split $h(y\cap z)$.
Since $Q$ was arbitrary it follows that $\{h(y\cap z):z\in\mathscr{G}\}$ is unsplittable, so we are done.

\hfill \qedref{thmmt}

Any $\mathfrak{r}_\kappa$-famliy of sets can be translated to a collection of functions in ${}^\kappa\kappa$ with a certain property related to splitting.
This translation between sets and functions will be useful.
We need the following definition:

\begin{definition}
\label{defbig} Big families of functions. \newline 
Let $\kappa$ be a regular cardinal and let $\mathcal{H}\subseteq{}^\kappa\kappa$. \newline 
The family $\mathcal{H}$ will be called big iff for every $\mathscr{F}\subseteq[\kappa]^\kappa$ and every $g\in{}^\kappa\kappa$ such that $\mathscr{F}$ contains the set $\{\beta\in\kappa:f(\beta)\leq g(\beta)\}$ whenever $f\in\mathcal{H}$ one can find an increasing $<\kappa$-to-one function $h\in{}^\kappa\kappa$ for which $\{h(y\cap z):y,z\in\mathscr{F}\}$ is an $\mathfrak{r}_\kappa$-family.
\end{definition}

For converting sets into functions we shall use a kind of projection.
Suppose that $y\in[\kappa]^\kappa$.
We define a function $p_y\in{}^\kappa\kappa$ by letting $p_y(\alpha) = \min(y\cap[\alpha,\kappa))$ for every $\alpha\in\kappa$.

\begin{claim}
\label{clm2big} Let $\kappa$ be a regular cardinal and let $\mathcal{R}$ be an $\mathfrak{r}_\kappa$ family. \newline 
The collection $\mathcal{H} = \{p_y:y\in\mathcal{R}\}$ is big.
\end{claim}

\par\noindent\emph{Proof}. \newline 
Fix any $\mathfrak{r}_\kappa$ family $\mathcal{R}$.
Assume toward contradiction that $\mathcal{H} = \mathcal{H}(\mathcal{R})$ is not big.
By definition, there are $g\in{}^\kappa\kappa$ and $\mathscr{F}\subseteq[\kappa]^\kappa$ such that $\mathscr{F}$ is upward-closed, $\{\beta\in\kappa:p_y(\beta)\leq g(\beta)\}\in\mathscr{F}$ whenever $y\in\mathcal{R}$ and for every $<\kappa$-to-one increasing $h\in{}^\kappa\kappa$ the family $\{h(y\cap z):y,z\in\mathscr{F}\}$ is not an $\mathfrak{r}_\kappa$-family.

By induction on $\alpha\in\kappa$ we define the interval $I_\alpha$ as follows. For $\alpha=0$ we simply take $I_0=[0,1)$. If $I_\alpha=[\gamma_\alpha,\gamma_{\alpha+1})$ has been defined then we choose $\gamma_{\alpha+2}\in\kappa$ such that $\forall\delta<\gamma_{\alpha+1}, g(\delta)<\gamma_{\alpha+2}$ and we let $I_{\alpha+1} = [\gamma_{\alpha+1},\gamma_{\alpha+2})$.
Finally, if $\alpha$ is a limit ordinal and $I_\beta=[\gamma_\beta,\gamma_{\beta+1})$ has been defined for every $\beta<\alpha$ then we let $\gamma_\alpha=\bigcup_{\beta<\alpha}\gamma_{\beta+1}$.
We choose $\gamma_{\alpha+1}\in\kappa$ such that $\forall\delta<\gamma_{\alpha}, g(\delta)<\gamma_{\alpha+1}$ and we let $I_\alpha=[\gamma_\alpha,\gamma_{\alpha+1})$.

Let $\Pi = \{I_\alpha:\alpha\in\kappa\}$.
Apply Theorem \ref{thmmt} to the triple $(\mathcal{R},\mathscr{F},\Pi)$ here standing for $(\mathscr{F},\mathscr{G},\Pi)$ there.
Among the three options given in Theorem \ref{thmmt}, $(a)$ and $(c)$ are excluded.
Firstly we show that $\mathcal{R}$ cannot be $\Pi$-scattered.
For this end, fix any nice $\Pi$-orbit $(D,E)$ and let $S=D$.
Notice that $E\subseteq(\kappa-S)$.
If $\mathcal{R}$ is $\Pi$-scattered then for every $y\in\mathcal{R}$ we have $|S\cap y|=|(\kappa-S)\cap y|=\kappa$.
This means that $S$ splits $\mathcal{R}$ which is impossible since $\mathcal{R}$ is an $\mathfrak{r}_\kappa$-family.
Secondly, our assumption toward contradiction (as unfolded in the first paragraph of the proof) says that $(c)$ of Theorem \ref{thmmt} fails.
We conclude, therefore, that $(b)$ of Theorem \ref{thmmt} holds, i.e. $\Pi$ is $\mathscr{F}$-scattered.

Let $(D,E)$ be a nice $\Pi$-orbit which exemplifies this fact.
Define $S = \bigcup\{I_\alpha\cup I_{\alpha+1}:I_\alpha\subseteq D\}$.
Similarly, let $T = \bigcup\{I_\alpha\cup I_{\alpha+1}:I_\alpha\subseteq E\}$.
Observe that $D\subseteq S$ and $E\subseteq T$ but still $S\cap T=\varnothing$ since there are no adjacent intervals in $D,E$.
Hence $T\subseteq(\kappa-S)$.

Fix any $y\in\mathcal{R}$ and let $A_y = \{\beta\in\kappa:p_y(\beta)\leq g(\beta)\}$, so $A_y\in\mathscr{F}$.
Since $\Pi$ is $\mathscr{F}$-scattered, $|A_y\cap D|=\kappa$.
But if $\beta\in A_y\cap D$ then $p_y(\beta)\in y\cap S$ since $I_{\alpha+1}\subseteq S$ whenever $I_\alpha\subseteq D$.
We conclude that $|y\cap S|=\kappa$.
Similarly, $|A_y\cap E|=\kappa$ and hence $|y\cap T|=\kappa$ which implies $|y\cap(\kappa-S)|=\kappa$.
It follows that $S$ splits all the elements of the $\mathfrak{r}_\kappa$-family $\mathcal{R}$, a contradiction.

\hfill \qedref{clm2big}

We need one last concept before proving the theorem below.

\begin{definition}
\label{deffindominating} Finite domination. \newline 
Let $\kappa$ be a regular cardinal and $\mathcal{A}\subseteq[\kappa]^\kappa$. \newline 
We say that $\mathcal{A}$ is finitely dominating iff for every $h\in{}^\kappa\kappa$ there is some finite collection $\{f_1,\ldots,f_n\}\subseteq\mathcal{A}$ such that $h\leq^*\max\{f_1,\ldots,f_n\}$.
\end{definition}

If $\mathcal{A} = \{f_\alpha:\alpha<\lambda\}\subseteq[\kappa]^\kappa$ and $\lambda<\mathfrak{d}_\kappa$ then $\mathcal{A}$ is not a dominating family and moreover it is not finitely dominating since $[\lambda]^{<\omega}=\lambda$.
We shall use this fact for proving the following theorem.

\begin{theorem}
\label{thmrukappa} If $\mathfrak{r}_\kappa<\mathfrak{d}_\kappa$ then $\mathfrak{r}_\kappa=\mathfrak{u}_\kappa$. \newline 
Consequently, $\mathfrak{r}_\kappa<\mathfrak{d}_\kappa$ implies $\cf(\mathfrak{r}_\kappa)>\omega$.
\end{theorem}

\par\noindent\emph{Proof}. \newline 
Let $\lambda=\mathfrak{r}_\kappa<\mathfrak{d}_\kappa$.
Choose an $\mathfrak{r}_\kappa$ family $\mathcal{R}$ of size $\lambda$.
We shall construct an ultrafilter $\mathscr{U}$ with a base of size $\lambda$, thus proving that $\mathfrak{u}_\kappa\leq\mathfrak{r}_\kappa$.
Since $\mathfrak{r}_\kappa\leq\mathfrak{u}_\kappa$ is always true we will be done.

Set $\mathcal{H} = \{p_y:y\in\mathcal{R}\}$.
From Claim \ref{clm2big} we infer that $\mathcal{H}$ is big.
Since $|\mathcal{H}|\leq\lambda<\mathfrak{d}_\kappa$ it is not finitely dominating.
Hence we can fix a function $g\in{}^\kappa\kappa$ which is not dominated by any finite number of functions from $\mathcal{H}$.
Define:
$$
\mathscr{B} = \{\{\beta\in\kappa:f(\beta)\leq g(\beta)\}:f\in\mathcal{H}\}.
$$
Let us point to some simple properties of $\mathscr{B}$.
First, this collection of sets has the strong finite intersection property.
This follows from the fact that $\mathcal{H}$ is not finitely dominating.
Second, the pair $(\mathscr{B},g)$ satisfies the assumptions in the definition of a big set with respect to $\mathcal{H}$, where $\mathscr{B}$ stands for $\mathscr{F}$ in the definition.
Consequently, there is an increasing $<\kappa$-to-one function $h\in{}^\kappa\kappa$ such that $\{h(y\cap z):y,z\in\mathscr{B}\}$ is an $\mathfrak{r}_\kappa$ family.
Finally, the cardinality of $\mathscr{B}$ is at most $\lambda$, since $|\mathcal{H}|\leq\lambda$.

Extend $\mathscr{B}$ to any ultrafilter $\mathscr{U}$ over $\kappa$ and notice that a base for $\mathscr{U}$ can be obtained from the elements of $\mathscr{B}$ by taking finite intersections. It follows that this base is of size at most $\lambda$ (and hence equals $\lambda$ since $\lambda=\mathfrak{r}_\kappa\leq\mathfrak{u}_\kappa$).
This observation concludes the proof.

\hfill \qedref{thmrukappa}

In the main result of the next section, $\mathfrak{r}_\lambda<\mathfrak{d}_\lambda$ is forced upon a singular cardinal $\lambda$.
We do not know whether $\mathfrak{r}_\lambda<\mathfrak{d}_\lambda$ implies $\mathfrak{r}_\lambda=\mathfrak{u}_\lambda$ when $\lambda>\cf(\lambda)$, though this is the typical case in the models of the next section.
This invites the following:

\begin{question}
\label{qru} Assume that $\lambda>\cf(\lambda)$. 
\begin{enumerate}
\item [$(\aleph)$] Is it consistent that $\mathfrak{r}_\lambda<\mathfrak{u}_\lambda$?
\item [$(\beth)$] Is it provable that $\mathfrak{r}_\lambda<\mathfrak{d}_\lambda$ implies $\mathfrak{r}_\lambda=\mathfrak{u}_\lambda$?
\end{enumerate}
\end{question}

Remark that the proof of $\mathfrak{r}_\kappa<\mathfrak{d}_\kappa$ implies $\mathfrak{r}_\kappa=\mathfrak{u}_\kappa$ when $\kappa$ is regular translates, for the most part, to a singular cardinal $\lambda$.
It seems, however, that the proof of Claim \ref{clm2big} is problematic in the case of a singular cardinal.
Specifically, the construction of the intervals $I_\alpha$ need not be a partition of $\lambda$ to $\lambda$-many intervals.

\newpage 

\section{At singular cardinals}

In this section we prove the consistency of $\mathfrak{u}_\lambda<\mathfrak{d}_\lambda$ where $\lambda>\cf(\lambda)$ is a strong limit cardinal.
In the definition of $\mathfrak{r}_\lambda$ or $\mathfrak{u}_\lambda$ there is no difference between the regular and the singular case.
But the definition of $\mathfrak{d}_\lambda$ requires some attention.
If $\lambda$ is singular then being of size less than $\lambda$ and being bounded in $\lambda$ are not the same statement.
We shall use the size version, as articulated in Definition \ref{defdkappa}.

Assume that $\mu>\cf(\mu)=\theta$.
The concept of $\mathfrak{d}_\mu$ relates to functions from $\mu$ into $\mu$, but the following useful lemma shows that one can deal with functions from $\mu$ into $\theta$.
Given $f,g\in{}^\mu\theta$ we shall say that $f<^* g$ iff $\{\beta\in\mu: f(\beta)\geq g(\beta)\}$ is of size less than $\mu$.
Define $\mathfrak{d}^*_\mu$ as the minimal cardinality of a domintaing subset of ${}^\mu\theta$ with respect to $<^*$.

\begin{lemma}
\label{lemcof} Assume that $\mu>\cf(\mu)=\theta$ and $\mathfrak{d}^*_\mu=\kappa$. \newline 
Then $\mathfrak{d}_\mu=\kappa$ as well.
\end{lemma}

\par\noindent\emph{Proof}. \newline 
As a first step we show that $\mathfrak{d}^*_\mu=\kappa$ implies $\mathfrak{d}_\mu\geq\kappa$.
Let $\mathscr{F}\subseteq{}^\mu\mu$ be of size $\tau<\kappa$.
We must show that $\mathscr{F}$ is not a dominating family in ${}^\mu\mu$.
Fix an increasing sequence of regular cardinals $(\lambda_i:i\in\theta)$ such that $\theta<\lambda_0$ and $\bigcup_{i\in\theta}\lambda_i=\mu$.
Let $h:\mu\rightarrow\theta$ be the associated interval mapping.
Explicitly, if $\alpha\in\mu$ then $h(\alpha)$ is the unique ordinal $i\in\theta$ such that $\lambda_i\leq\alpha<\lambda_{i+1}$.

For every $f\in\mathscr{F}$ we define $f'\in{}^\mu\theta$ by $f'(\beta)=h(f(\beta))$ whenever $\beta\in\mu$.
Let $\mathscr{F}' = \{f':f\in\mathscr{F}\}$, so $|\mathscr{F}'|\leq\tau<\kappa$.
Since $\mathfrak{d}^*_\mu\geq\kappa$ there exists $g'\in{}^\mu\theta$ such that $\forall f'\in\mathscr{F}', \neg(g'\leq^* f')$.
Define $g\in{}^\mu\mu$ by letting $g(\beta)=\lambda_{i+1}\Leftrightarrow g'(\beta)=i$.
We claim that $g$ is not dominated by $\mathscr{F}$.

For proving this claim fix an element $f\in\mathscr{F}$ and recall that $\neg(g'\leq^* f')$.
Hence there exists a set $A\in[\mu]^\mu$ such that $\beta\in A\Rightarrow g'(\beta)>f'(\beta)$.
By the definition of $f'$ we may write $\beta\in A\Rightarrow g'(\beta)>h(f(\beta))$.
This means that $g(\beta) = \lambda_{h(f(\beta))+1}>f(\beta)$ whenever $\beta\in A$, so $\neg(g\leq^* f)$ and the first direction is proved.

For the opposite direction assume that $\mathfrak{d}_\mu=\kappa$, aiming to prove that $\mathfrak{d}^*_\mu\geq\kappa$.
If we succeed then assuming $\mathfrak{d}^*_\mu=\kappa$ one concludes that $\mathfrak{d}_\mu>\kappa$ is impossible, so $\mathfrak{d}^*_\mu=\kappa$ implies $\mathfrak{d}_\mu\leq\kappa$ as required for this direction.
Let $\mathscr{G}\subseteq{}^\mu\theta$ be a family of elements of ${}^\mu\theta$ such that $|\mathscr{G}|=\tau<\kappa$.
We shall see that $\mathscr{G}$ is not dominating.

For each $g\in\mathscr{G}$ we define $f_g\in{}^\mu\mu$ by letting $f_g(\beta) = \lambda_{h(g(\beta))}$.
Set $\mathscr{F} = \{f_g:g\in\mathscr{G}\}$.
Since $|\mathscr{F}|\leq\tau<\kappa$, we can choose $g^{\rm up}\in{}^\mu\mu$ such that $\forall f\in\mathscr{F}, \neg(g^{\rm up}\leq^* f)$.
We define a function $g\in{}^\mu\theta$ as follows.
If $\beta\in\mu$ then there is a unique ordinal $i\in\theta$ for which $\lambda_i\leq g^{\rm up}(\beta)<\lambda_{i+1}$, and we let $g(\beta)=i+1$.
We claim that $g$ is not dominated by $\mathscr{G}$.

By way of contradiction assume that $g_0\in\mathscr{G}$ and $g\leq^* g_0$.
Let $f_0$ be $f_{g_0}$, and recall that $A = \{\beta\in\mu: f_0(\beta)<g^{\rm up}(\beta)\}$ is of size $\mu$.
Let $B = \{\gamma\in\mu: g(\gamma)>g_0(\gamma)\}$ and notice that $|B|<\mu$.
Fix an ordinal $\gamma_0\in\mu$ such that $B\subseteq\gamma_0$.
If $\beta\in A-\gamma_0$ then $f_0(\beta)<g^{\rm up}(\beta)$ (since $\beta\in A$), which means that $\lambda_{h(g_0(\beta))}<\lambda_{h(g(\beta))}$ (by the definition of these functions), hence $g_0(\beta)<g(\beta)$.
This is impossible, however, since $\beta\geq\gamma_0$, so we are done.

\hfill \qedref{lemcof}

The first observation below is that $\mathfrak{d}_\lambda$ behaves nicely at singular cardinals.

\begin{claim}
\label{clmdgeqlambda} Let $\lambda$ be a singular cardinal. \newline 
Then $\mathfrak{d}_\lambda>\lambda$, and moreover $\cf(\mathfrak{d}_\lambda)>\lambda$.
\end{claim}

\par\noindent\emph{Proof}. \newline 
We spell out the easy argument for $\mathfrak{d}_\lambda>\lambda$, and it will also follow of course from the stronger statement $\cf(\mathfrak{d}_\lambda)>\lambda$.
Suppose that $\{f_\alpha:\alpha\in\lambda\}\subseteq{}^\lambda\lambda$.
We wish to describe a function $h:\lambda\rightarrow\lambda$ which is not dominated by any $f_\alpha$.
For this end, decompose $\lambda$ into $\lambda$-many disjoint sets $(S_\alpha:\alpha\in\lambda)$, each of which of size $\lambda$.
For every $\beta\in\lambda$ let $h(\beta)=f_\alpha(\beta)+1$ iff $\alpha$ is the unique ordinal for which $\beta\in S_\alpha$.
For every $\alpha\in\lambda$ one can see that $\beta\in S_\alpha\Rightarrow f_\alpha(\beta)<h(\beta)$, so $h$ is as required.

Assume now that $\mathfrak{d}_\lambda=\chi$, and assume toward contradiction that $\theta=\cf(\chi)<\lambda$.
Let $\kappa = \cf(\lambda)$.
By Lemma \ref{lemcof} we may concentrate on ${}^\lambda\kappa$, so fix a family $\mathcal{F}\subseteq{}^\lambda\kappa$ which exemplifies $\mathfrak{d}_\lambda=\chi$.
Choose a sequence of disjoint sets $(\mathcal{F}_i:i\in\theta)$ such that $\mathcal{F}=\bigcup_{i\in\theta}\mathcal{F}_i$ and $|\mathcal{F}_i|<\chi$ for every $i\in\theta$.
Decompose $\mu$ into $(A_i:i\in\theta)$, each $A_i$ is of size $\mu$, and fix a bijection $g_i:\mu\rightarrow A_i$ for every $i\in\theta$.

Now for each $i\in\theta$ let $\mathcal{G}_i = \{f\circ g_i:f\in\mathcal{F}_i\}$.
Observe that $|\mathcal{G}_i|<\chi$ for every $i\in\theta$.
For every $i\in\theta$ we have $|\mathcal{F}_i|<\chi$, and hence one can choose a function $f_i\in{}^\lambda\kappa$ such that $\forall f\in\mathcal{F}_i, \neg(f_i\leq^* f)$.
Much as in the first part of the proof, we can describe our non-dominated function $g$ by defining its values separately over each $A_i$.
For every $i\in\theta$ and for every $\alpha\in A_i$ let $g(\alpha) = f_i(g_i^{-1}(\alpha))$.
Notice that $g\in{}^\lambda\kappa$ is well defined.
We claim that $g$ is not dominated by $\mathcal{F}$, namely $\forall f\in\mathcal{F}, \neg(g\leq^* f)$.

For proving this statement, fix an element $f\in\mathcal{F}$.
Let $i\in\theta$ be the unique ordinal for which $f\in\mathcal{F}_i$.
Denote the set $\{\beta\in\mu:f\circ g_i(\beta)<f_i(\beta)\}$ by $B$, so $|B|=\mu$.
Let $C = \{g_i(\beta):\beta\in B$.
Observe that $C\subseteq A_i$ and $|C|=\mu$ since $g_i$ is one-to-one.
We claim that $f(\gamma)<g(\gamma)$ whenever $\gamma\in C$.
Indeed, fix any element $\gamma\in C$, and let $\beta\in B$ be the unique ordinal for which $\gamma=g_i(\beta)$.
By the definition of $B$ it follows that $f(\gamma)=f(g_i(\beta))<f_i(\beta)$.
Concomitantly, $g(\gamma)=f_i(g_i^{-1}(\gamma))=f_i(\beta)$, and hence $f(\gamma)<g(\gamma)$.
Since $|C|=\mu$ we conclude that $\neg(g\leq^* f)$.
But $f$ was arbitrary, so we are done.

\hfill \qedref{clmdgeqlambda}

For proving the main result of this section we need some controll over the true cofinality of some sequences of regular cardinals.
Aiming to show that $\mathfrak{u}_\lambda<\mathfrak{d}_\lambda$ is consistent, one sequence or regular cardinals will give a (small) upper bound for $\mathfrak{u}_\lambda$ and the other sequence will give a (large) lower bound for $\mathfrak{d}_\lambda$.
More precisely, for every characteristic we must force the value of the true cofinality for both the sequence of regular cardinals and the sequence of their successors.
The forcing machinery for this end comes from \cite{MR2987137}, and we shall use the following version:

\begin{theorem}
\label{thmtcf} Let $\lambda$ be a supercompact cardinal. \newline 
Then one can force the following statements:
\begin{enumerate}
\item [(a)] $\lambda>\cf(\lambda)=\theta$.
\item [(b)] $\lambda<\kappa=\cf(\kappa)\leq 2^\lambda$.
\item [(c)] $(\lambda_i:i<\theta)$ is an increasing sequence of strongly inaccessible cardinals, $\theta<\lambda_0$ and $\lambda=\bigcup_{i\in\theta}\lambda_i$.
\item [(d)] $2^{\lambda_i}=\lambda_i^+$ for every $i\in\theta$.
\item [(e)] ${\rm tcf}(\prod_{i\in\theta}\lambda_i,J^{\rm bd}_\theta)=\kappa$.
\item [(f)] ${\rm tcf}(\prod_{i\in\theta}\lambda_i^+,J^{\rm bd}_\theta)=\kappa$.
\end{enumerate}
Moreover, $2^\lambda$ can be arbitrarily large, and $\kappa$ can be an arbitrarily large regular cardinal provided that $\kappa\leq 2^\lambda$.
\end{theorem}

\hfill \qedref{thmtcf}

We show now how to take care of the dominating number:

\begin{theorem}
\label{thmbigd} Assume that $\lambda$ is supercompact. \newline 
Then one can force $\lambda$ to be a strong limit singular cardinal, $2^\lambda$ is arbitrarily large, $\kappa=\cf(\kappa)\leq 2^\lambda$ is arbitrarily large above $\lambda$ and $\mathfrak{d}_\lambda\geq\kappa$.
\end{theorem}

\par\noindent\emph{Proof}. \newline 
Apply Theorem \ref{thmtcf} to obtain $\theta=\cf(\lambda)<\lambda$, and let $(\lambda_i:i<\theta)$ be an increasing sequence of strongly inaccessible cardinals such that $\theta<\lambda_0$ and $\lambda=\bigcup_{i\in\theta}\lambda_i$ as guaranteed there.
It means that $2^{\lambda_i}=\lambda_i^+$ for every $i\in\theta$, with both ${\rm tcf}(\prod_{i\in\theta}\lambda_i,J^{\rm bd}_\theta)=\kappa$ and 
${\rm tcf}(\prod_{i\in\theta}\lambda_i^+,J^{\rm bd}_\theta)=\kappa$.
Denote $J^{\rm bd}_\theta$ by $J$.
Fix a sequence $(f_\alpha:\alpha\in\kappa)$ of functions in the product $\prod_{i\in\theta}\lambda_i$ which exemplifies ${\rm tcf}(\prod_{i\in\theta}\lambda_i,J)=\kappa$.
Fix also a sequence $(g_\alpha:\alpha\in\kappa)$ of functions in the product $\prod_{i\in\theta}\lambda_i^+$ which exemplifies ${\rm tcf}(\prod_{i\in\theta}\lambda_i^+,J)=\kappa$.

For every $i\in\theta$ enumerate the elements of ${}^{\lambda_i}\theta$ by $\mathscr{F}_i = \{g^i_\alpha:\alpha\in\lambda_i^+\}$.
Likewise, for each $i\in\theta$ fix a sequence of mappings $(h^i_\alpha:\alpha\in\lambda_i^+)$ so that every $h^i_\alpha$ is a one-to-one mapping from $\alpha$ into $\lambda_i$.
For every $\alpha\in\kappa$ and every $i\in\theta$, define:
$$
w_{\alpha i} = \{\beta\in\lambda_i^+: \beta<g_\alpha(i) \wedge h^i_{g_\alpha(i)}(\beta)<f_\alpha(i)\}.
$$
Observe that $w_{\alpha i}\in[\lambda_i^+]^{<\lambda_i}$ for every $\alpha\in\kappa, i\in\theta$.
Indeed, $w_{\alpha i}\subseteq\lambda_i^+$ by definition, and its cardinality is bounded by $|f_\alpha(i)|<\lambda_i$.
The focal property of the sets $w_{\alpha i}$ is that if $(u_i:i\in\theta)$ satisfies $u_i\in[\lambda_i^+]^{<\lambda_i}$ for every $i\in\theta$ then for some $\alpha_1\in\kappa$, if $\alpha\in(\alpha_1,\kappa)$ then there is $i(\alpha)\in\theta$ such that for every $i\in(i(\alpha),\theta)$ it is true that $u_i\subseteq w_{\alpha i}$.

For proving this property, fix a sequence $(u_i:i\in\theta)$ such that $u_i\in[\lambda_i^+]^{<\lambda_i}$ for every $i\in\theta$.
For each $i\in\theta$ let $g(i)=\sup(u_i)$, so $g\in\prod_{i\in\theta}\lambda_i^+$.
By the choice of $(g_\alpha:\alpha\in\kappa)$ there is an ordinal $\alpha_0\in\kappa$ such that if $\alpha\in(\alpha_0,\kappa)$ then $g\leq^* g_\alpha$.
Now for every $i\in\theta$ let $v_i = \{h^i_{g_\alpha(i)}(\beta):\beta\in u_i\}$. The cardinality of $v_i$ is less than $\lambda_i$ and hence it is bounded in $\lambda_i$.
Define $f(i)=\sup(v_i)$ for every $i\in\theta$, so $f\in\prod_{i\in\theta}\lambda_i$.
Fix an ordinal $\alpha_1\in[\alpha_0,\kappa)$ such that if $\alpha\in(\alpha_1,\kappa)$ then $f<^* f_\alpha$.
Now if $\alpha\in(\alpha_1,\kappa)$ then a sufficiently large $i_1\in\theta$ such that $i\in(i_1,\theta) \wedge \beta\in u_i \Rightarrow [\beta<g_\alpha(i)\wedge h^i_{g_\alpha(i)}(\beta)<f_\alpha(i)]$, which amounts to $u_i\subseteq w_{\alpha i}$, so we are done.

Back to the main argument, for every ordinal $\alpha\in\kappa$ we choose a function $h_\alpha:\lambda\rightarrow\theta$ with the following property:
\begin{center}
If $\beta\in w_{\alpha i}, j<\theta$ and $|g^i_\beta(j)|=\lambda_i$ \\
then $|h_\alpha^{-1}((j,\theta))\cap g^i_\beta(j)|=\lambda_i$.
\end{center}
For choosing these functions notice that one can take care of a single $j\in\theta$ and a single ordinal $\beta$, since $|g^i_\beta(j)|=\lambda_i$.
Now for each $\alpha\in\kappa$ there are only $\theta\times|w_{\alpha i}|<\lambda_i$ many pairs of the form $(j,\beta)$ to take care of, so the choice of these functions is possible.

Let $g$ be any function from $\lambda$ into $\theta$.
If there are $\theta$-many ordinals $i$ for which $g\upharpoonright\lambda_i\in\{g^i_\beta:\beta\in w_{\alpha i}\}$ then $\neg(h_\alpha\leq^* g)$ by the above property of the $h_\alpha$s.
Given a function $g\in{}^\lambda\theta$, let $\beta_i$ be the unique ordinal so that $g\upharpoonright\lambda_i=g^i_{\beta_i}$, and let $u^i_g=\{\beta_i\}$ for every $i\in\theta$.
By the focal property of the $w_{\alpha i}$s, for some $\alpha_1\in\kappa$, if $\alpha\in(\alpha_1,\kappa)$ then there is $i(\alpha)\in\theta$ such that for every $i\in(i(\alpha),\theta)$ it is true that $u^i_g\subseteq w_{\alpha i}$.
At each ordinal $\alpha\in(\alpha_1,\kappa)$ we see, therefore, that $\neg(h_\alpha\leq^* g)$.

Assume, now, that $\mathscr{G}\subseteq{}^\lambda\theta$ and $|\mathscr{G}|=\tau<\kappa$.
Enumerate the elements of $\mathscr{G}$ by $\{g_\delta:\delta\in\tau\}$.
For every $\delta\in\tau$ choose $\alpha_\delta\in\kappa$ such that if $\alpha\in(\alpha_\delta,\kappa)$ then $\neg(h_\alpha\leq^* g_\delta)$.
Let $\alpha = \bigcup_{\delta\in\tau}\alpha_\delta$, so $\alpha<\kappa$.
It follows that $\forall g\in\mathscr{G}, \neg(h_\alpha\leq^* g)$, so $\mathscr{G}$ is not a dominating family.
Since $\mathscr{G}$ was an arbitrary family of less than $\kappa$ many functions in ${}^\lambda\theta$ we conclude that $\mathfrak{d}_\lambda\geq\kappa$, thus accomplishing the proof.

\hfill \qedref{thmbigd}

It has been proved in \cite{MR2992547} that one can force $\mathfrak{u}_\lambda=\lambda^+$ and $2^\lambda$ is arbitrarily large for some singular cardinal $\lambda$.
The proof is similar to the proof in this section, but somehow from the opposite direction. Namely, under similar assumptions about the true cofinalities of sequences of regular cardinals and their successors, one can show that $\mathfrak{u}_\lambda\leq\kappa$ where $\kappa$ realizes the true cofinalities.
This is possible, in particular, for $\kappa=\lambda^+$.
It remains to merge the result about $\mathfrak{u}_\lambda$ with the present result about $\mathfrak{d}_\lambda$.
For this end, one needs two sequences of regular cardinals, with different values of true cofinalities.
The ability to force this situation appears in \cite{1040} and \cite{MR3201820}.

\begin{theorem}
\label{thmuandd} Assume the existence of a supercompact cardinal. \newline 
Then one can force a singular cardinal $\lambda>\cf(\lambda)=\theta$ such that $\lambda<\mathfrak{u}_\lambda\leq\kappa_0<\kappa_1\leq\mathfrak{d}_\lambda\leq 2^\lambda$, and the gap between $\mathfrak{u}_\lambda$ and $\mathfrak{d}_\lambda$ can be arbitrarily large.
\end{theorem}

\par\noindent\emph{Proof}. \newline 
For $\lambda>\cf(\lambda)=\theta$ and $\kappa=\cf(\kappa)\in(\lambda,2^\lambda]$ we say that $(\lambda_i:i\in\theta)$ is $\kappa$-qualified if it is an increasing sequence of regular cardinals such that $\theta<\lambda_0, \lambda=\bigcup_{i\in\theta}\lambda_i, 2^{\lambda_i}=\lambda_i^+$ for every $i\in\theta$ and ${\rm tcf}(\prod_{i\in\theta}\lambda_i,J) = {\rm tcf}(\prod_{i\in\theta}\lambda_i^+,J)=\kappa$ (where $J$ is usually the ideal of bounded subsets of $\theta$, but not necessarily).

Let $\lambda$ be supercompact in the ground model.
Let $\mathbb{P}$ be a forcing notion which makes $\lambda>\cf(\lambda)=\theta$ and forces the following statements:
\begin{itemize}
\item There exists a sequence $(\lambda_i^0:i\in\theta)$ of measurable cardinals which is $\kappa_0$-qualified.
\item There exists a sequence $(\lambda_i^1:i\in\theta)$ of strongly inaccessible cardinals which is $\kappa_1$-qualified.
\end{itemize}
From Theorem 1.4 of \cite{MR2992547} we infer that $\mathfrak{r}_\lambda\leq\mathfrak{u}_\lambda\leq\kappa_0$.
From Theorem \ref{thmbigd} we infer that $\mathfrak{d}_\lambda\geq\kappa_1$.
since $\lambda<\mathfrak{r}_\lambda$ and $\mathfrak{d}_\lambda\leq 2^\lambda$ are always true, we are done.

\hfill \qedref{thmuandd}

The ability to force two such qualified sequences is a version of Theorem \ref{thmtcf}, see \cite{MR3201820}.
Similar statements appear in \cite{1040}, in a slightly different way.
In that paper, many regular cardinals above the singular cardinal $\lambda$ are realized as true cofinalities of some sequence of measurable cardinals below $\lambda$.
However, if ${\rm tcf}(\prod_{i\in\theta}\lambda_i,J)=\kappa$ then ${\rm tcf}(\prod_{i\in\theta}\lambda_i^+,J)=\kappa^+$.
For our purpose it makes no difference, since we will get $\mathfrak{u}_\lambda\leq\kappa_0^+$ and $\mathfrak{d}_\lambda\geq\kappa_1^+$, still there will be a gap between these characteristics.
But there is an important difference, since in \cite{1040} the GCH is kept below $\lambda$ while $2^\lambda>\lambda^+$, and being strong limit this means that we must work with a singular cardinal with countable cofinality.
On the other hand, this theorem is very flexible in the sense that an infinite set of targets obtains by qualified sequences.
This is the background behind the following statement.

\begin{corollary}
\label{cordsing} It is consistent that $\lambda>\cf(\lambda), \mathfrak{u}_\lambda<\mathfrak{d}_\lambda$ and $\mathfrak{d}_\lambda$ is singular.
\end{corollary}

\par\noindent\emph{Proof}. \newline 
Let $\lambda$ be supercompact in the ground model.
We set $\kappa=\lambda^+$ and we choose an increasing sequence of regular cardinals $(\kappa_j:j\in\partial)$ where $\partial\geq\lambda^+, \cf(2^\lambda)=\partial$ in the generic extension and $\bigcup_{j\in\partial}\kappa_j=2^\lambda$.
We shall say that the sequence $(\lambda_n:n\in\omega)$ is $(\theta,\theta^+)$-qualified iff $2^{\lambda_n}=\lambda_n^+$ for every $n\in\omega$, $\lambda=\bigcup_{n\in\omega}\lambda_n, J\supseteq J^{\rm bd}_\omega$ and $\theta={\rm tcf}(\prod_{n\in\omega}\lambda_n,J)$ while $\theta^+={\rm tcf}(\prod_{n\in\omega}\lambda_n^+,J)$.

We repeat the argument of the previous theorem, but now we employ Theorem 7 from \cite{1040}.
In particular, $\lambda>\cf(\lambda)=\omega$ and a $(\kappa_j,\kappa_j^+)$-qualified sequence of measurable cardinals below $\lambda$ is forced for every $j\in\partial$ as well as a sequence for $(\kappa,\kappa^+)$.
It follows from Theorem \ref{thmbigd} that $\mathfrak{d}_\lambda\geq\kappa_j$ for every $j\in\partial$.
Hence $\mathfrak{d}_\lambda=2^\lambda$, so $\mathfrak{u}_\lambda\leq\kappa^+<\mathfrak{d}_\lambda$ and $\mathfrak{d}_\lambda$ is singular as required.

\hfill \qedref{cordsing}

\newpage 

\bibliographystyle{amsplain}
\bibliography{arlist}

\end{document}